\documentclass{article}\usepackage{graphicx}\usepackage{color} \usepackage{amsmath}\usepackage{amsfonts}\usepackage{newce
nt}\usepackage{graphicx}\begin{document}\large
\definecolor{x}{RGB}{10,10,10}
\definecolor{y}{RGB}{250,245,241}
\def \m {\vspace{12pt}

\noindent}
\def \R {{\mathbb R}}
\def \C {{\mathbb C }}
\def \N {{\mathbb N}}
\def \Q {{\mathbb Q}}
\def \Z {{\mathbb Z}}
\color{x}
\pagecolor{y}
 {\small  
\pagenumbering{roman}
\centerline{Arc lifting  for the Nash manifold}
\m 
\m
\m
\m
\m
\m
\m \hfill  $\begin{matrix}
 \hbox{\it John Atwell Moody}\hfill \cr 
\hbox{\it Warwick Mathematics Institute}\hfill \cr\hbox{\it  March, 2011}\hfill

\end{matrix}$
\m

}

\vfill\eject\noindent
{\bf Preface}
\m 
Our understanding of nature doesn't do anything
like encompass nature. 

\m
In  `The olden times' Ian Hislop [1] 
claims that Tolkien's romanticization of nature is misleading,  that the home town
which is implicitly the real setting of his fantasies
had already been industrialized  before he lived there.  
What Tolkien writes about actually does make sense. Not when you go to his
town as an investigator, it is when you think
 of stories over the generations. 
\m
Two of my grandparents
lived
in a  mill-town 
 on a river, with its `mill-pond' upstream from the water-wheel,
at the edge of
what they called the `never-ending forest.'
You could see Northern lights at night, and 
the  people living there, some of them immigrants from Norway and other northern places,
 knew,  loved
and understood the trees,  mosses and types of fungus shapes
that grow on healthy trees. Also the types of `grandaddy long-legs' spiders,
salamanders, and all types of living things. These live on
one side of the never-ending forest,  the side
which will never be completely cut down. As the boundary
moved northwards, the goodness of the never-ending forest seemed like it was being absorbed
into their lives.
\vfill\eject\noindent
{\bf Context of the Nash-Semple question}
\m
Nash's question in $n$ dimensions is this:  each time the indeterminacy of $n$ forms
is resolved, the regular locus of the resolution is an open subset
which extends the regular locus of a singular manifold, yielding 
a topological space with an open cover by
an  ascending series of open  manifolds.
\m
As a question in the classical
topology, it is
whether a continuous path starting in the original manifold,
and whose whole image there is compact, also
has an   endpoint in one of the open
subsets whose union is the Nash manifold, or, whether
it can actually go on forever.
\m
Resolving the indeterminacy of $n$-forms starts with thinking about
smoothly parametrized paths in space. The image of a path retains
some aspects of an implicit parametrization. An image in projective space of
a compact and connected Riemann surface under a holomorphic map,
if it is more than one point, is an algebraic curve and the Riemann
surface can be rediscovered up to an automorphism by normalizing the
algebraic curve. Normalization has a very algebraic definition,
a local section of the structure sheaf of the normalization is
nothing but a local section of the sheaf of endomorphisms of a torsion-free coherent rank one
sheaf over the original scheme (in this case the algebraic curve),
and also, as was proved by Nobile, if we resolve the indeterminacy
of the Gaussian map, sending each regular point of the curve to
its projective tangent line, and then repeat using the suitably
modified ambient manifold, we reach the normalization in finitely
many steps. And within the rational function field of an irreducible
algebraic variety (in this case our algebraic curve),
each discrete valuation is an `ideal point'
corresponding to a prime divisor of the normalization.
\vfill\eject\noindent
At each stage of Nash's process, in arbitrary dimensions, the
smooth locus is a smooth quasi-projective variety. 
It is worth thinking about why such a question should matter to anyone.

\m
Mandelbrot insisted that things in mathematics
which are nothing like manifolds may equally represent what is really
there, such as a coastline [2].  

\m
In hard science, students trying to learn quantum theory are told
that there is a chain of prerequisite definitions.
The ground state component of Hydrogen, initially imagined
as belonging to the real one-dimensional vector-space spanned
by the function $e^{-r}$ in suitable units of measurement,
needed to be re-imagined as belonging instead to the two-dimensional representation
of $SU_2$ which only exists over the complex numbers. If someone
said, ``It's one of Hamilton's quaternions, then,'' it would lead
to incorrectly believing in things like entanglement.
\m
After Mumford left Mathematics to think
collaboratively about
digital technology, he had written that he had
loved making mappings,  had
thought of it as a way to explore real things [3].
\m
When someone says that  all of reality is a manifold (or a singular
manifold or a scheme, or a stack, or an etale site) they are   failing
to see the tentative and provisional nature which mathematical definitions
must have.
\m

\m
Witten 
writes in  in  {\it physics and geometry} ``Spacetime is a pseudo-Riemannian manifold M, endowed with a metric tensor'' [4].
\m
Subsequent thinking  of Mori and Reid,  restricted specifically to
quasiprojective complex varieties,
says that   a smooth complex variety has
a  minimal model with {\it canonical} singularities,
which is locally a canonical model and more natural than a smooth
model [5].
\vfill\eject\noindent
Our thinking is doomed if we try to first make an axiomatic map of possible
worlds, and then find `nature' as an element on our map.
The tragedy is that we 
have no other choice
 if we are going to 
choose and enact, or decline to enact,  policies. We  must
rely on a concept of truth.
We must make definitions. Clarity of the definitions
implies axioms and consistency, and a mathematical map.
\m 
The necessary  consequence  is that people have to rely on the mathematics which
is most broadly accepted, like facts about the real number line. 
\m I am going not going to include  a political 
discourse, or describe the true depth of the tragedy, which can only be understood
through evolutionary psychology, and through seeing in every example
how things like the internet cannot couch coherent thought;
but only to say that focussing on global warming
does incorrectly  make it seem right to install agriculture and housing
on ancient wilderness nature, as long as we  compensate by  planting 
massive numbers of seedling
trees.  
\m
Scientific articles which  appear nowadays
say that  worrying about the accumulations in nature of
mutations chosen by people correlates with lacking scientific
understanding.
That correct scientific understanding, in turn, correlates statistically with
believing
that each next generation of corporate, government, and academic scientists
in every country know how to   mutate nature based on their
best goals, strategies and  biological models [6].
\m And yet, just as in quantum theory there was a dawning
realization that there is no `algebra' like a number line which
makes up an underlying universe, and that nor do Hydrogen waveforms really
 have a
quaternionic component, there is also in biology the notion that
there aren't primarily `genes' and `reading frames' which, in the
analogously trivial  way, would form any sensible  ambient structure pre-existing our
attempt to imagine it.
\vfill\eject\noindent
The best biologists know how to give only a hint of an idea, 
 of the type of saying, things you never thought could
be related to each other,  are related in surprising
and direct ways, but infinitely complex and previously unknown ways. 
That, if you have made a permanent mutation
without understanding {\it this}, you would have caused a tragedy.
\m It doesn't mean that we can call a halt to the production,
and enact legislation about {\it this}. It means that our understanding
is and always will be nothing but provisional. 

\m The effective refutation
on that front arises from a surprising place: not from
 advanced mathematical theory, but from  the GCSE curriculum for children.
Children are being encouraged  to wonder:  while there
may be  `the gene to increase fatty acid in pork,' whether there is not going to
be ever found `the gene which make an orchid resemble a bee,' or `the gene which
controls symbiosis.' 
\m
People who live in cities, or in established
countries, an increasing majority, who are said to have been `lifted' out of poverty, have a trivialized and almost mathematical
understanding of nature. They are like competent gardeners, 
and  things are dismissed as not important in the way that gardeners
dismiss weeds, or call grasshoppers `locusts.'  
They are the ones who will outnumber everyone and
prevail, and they will in all sincerity take a medical point of view.
Like how the puffins on the Isle of Lundy are `protected' until Lundy
is a puffin farm. 
\m
Chasing the effect of repeatedly resolving the
indeterminacy of $n$-forms to see where it would lead
could also just be doing also the wrong thing.
Mumford  is right to think that
just because we can make a map, doesn't in any way imply that
we can rely on the map to tell us what is really there.
I worry -- selfishly maybe -- that 
if mathematicians stop making {\it new} maps, then the last
ones which were made will become trusted as if they were
the only ones there, and thus had been a {\it necessary part}
of existence.
\m Whereas also, people who live in new, or not yet `developed' countries,
like America was not long ago, while they have a very deep
understanding of nature, one which can never be put into words,
and while they are horrified by anyone who thinks of development
as anything but an irreversible loss, yet, they see nature as
abundant, almost infinite, and forgiving. 
\m
Nash was an American,
and if it is true, as people say, that he wanted to construct
a limiting object out of resolving indeterminacy of $n$-forms,
then it would be wanting to `follow the yellow brick road.'
That it doesn't matter where it will go, but it will be to
a new and wonderful place which will make us happy. It doesn't matter
what our actions are, the results will be consistent and wonderful.
It is like the remaining people who still have `frontier' optimism,
in Alaska, or, I should say, now, rather,  in Northern Alaska.

\m
Subsequent to Nash's question, Hironaka did  answer
the question of resolving singularities (algebraically in characteristic
zero) [7]. His proof was based on reducing local multiplicities iteratively
and the particular resolution depended on a choice of embedding
into a smooth manifold.

\m
Perez and  Teissier  [8]
found an earlier reference [9] to Nash's question;  most references  say it is a question
of Nash, predating Hironaka's successful proof, with Nash putting
nothing in writing. Notable mathematicians have at times
informally made announcements of papers which (like mine about arc lifting)
have never  materialized.
\m
M. Spivakovsky has said, during a conversation in Grenoble [10], that
his and Hironaka's theorem can be understood as
a modification of Nash's question in which the surface is normalized
at each step.
\vfill\eject\noindent
The same day, S. Abyankhar [11] answered a question about finding a  local singularity
unaffected (up to isomorphism and further localization) by a Nash transform without saying anything, but,
to my best recollection of it,  acting like
trying to snatch a mosquito out of the air where there was no mosquito.
To a later question, what are the conditions for a Nash transform
to yield a normal manifold, he waited a significant time,  and only said,
it is a deep question [12].
\m
Gonzales-Sprinberg [13],  and later Atanasov,  Lopez, Perry,  Proudfoot and Thaddeus  have anayzed toric case [14].
\m
I'll also reference some things which I have not finished reading about,
Encinas and Villamayer [15] , Wlodarczyk [16], Kollar [17],  McQuillan [18],
and  Abramovich, Temkin and Wlodarczyk [19]
 have removed
dependence on embedding, though the later work needs stacks.
\m
Referees
while remaining anonymous have excused long delays in replying to
my submissions with saying they had been working themselves on Nash's
question, and in their hints to me have tended to introduce a differential geometry
feeling to the subject.
\vfill\eject\noindent
{\bf Primary Decomposition}
\m
One  way of understanding  cycles
is to view the structure sheaf of a closed subscheme  (or  analyic
subspace)
as a coherent sheaf on the ambient variety.
Let's restrict to the case of schemes;
in generalzing to the
analytic setting we'll need to remember  that a subsheaf of a coherent sheaf
need not be coherent  (while the kernel of
a map of coherent sheaves always is).
\m
A coherent sheaf ${\cal D}$ on a quasiprojective scheme is `coprimary' (generalizing the definition given
in `Algebre Local, Multiplicities' [20])
if there is a prime ideal sheaf ${\cal P}$ so that
if there is any prime ideal sheaf ${\cal Q}$ and
invertible sheaf ${\cal L}$ and embedding
${\cal L}\otimes {\cal O}_M/{\cal Q} \to {\cal C}$
then ${\cal P}={\cal Q}.$
\m
In any filtration
of such a  ${\cal C}$ with successive quotients of type
${\cal L}\otimes {\cal O}_M/{\cal Q}, $ the prime ideal sheaf 
 ${\cal P}$ is the unique minimal element in the set of ${\cal Q}.$  
The occurrence of ${\cal L}$ in the
definition is needed so that we can be sure that ${\cal C}$
does actually have such a filtration.
\m
Within any  coherent sheaf ${\cal C}$ we may choose a maximal coherent subsheaf
which is a  sum (necessarily a direct sum) of coprimary sheaves, and
I would guess that if we tensor  this with the cartesian product of the local rings
corresponding to the corresponding primes, the tensor product of
the subsheaf
should be isomorphic with the maximal Artinian subsheaf of the tensor
product of all of ${\cal C}$
\m
If there are no inclusions among these primes, we can use
instead the semi-localization, and the maximal artinian
subsheaf is the whole of the semilocalization of ${\cal C}.$
\vfill\eject\noindent
We also can use the completions instead of localizations, and
 the maximal Artinian module of the tensor product of our coherent
cartesian product of the localizations is already  a complete
module for the cartesian product of the completions.
It is complete for each cartesian factor a finite power of the radical
acts by zero.
\m Our maximal direct sum of coprimary subsheaves of ${\cal C}$ ought to be
{\it essential} in the sense that it intersects any nontrivial
subsheaf. \m
In this sense, subsheaves should be detected by discrete valuations.
The intersection of a subsheaf with our chosen ${\cal P}$-coprimary
subsheaf has a filtration with successive quotients
only of the type ${\cal L}\otimes {\cal O}_M/{\cal Q}$ for ${\cal Q}$
containing ${\cal P}$, and 
and we may count the occurrences of ${\cal P}$  directly, or just
tensor the intersection with the maximal coprimary subsheaf
with the finite cartesian product of localizations, and count the
irreducible composition factors of the resulting Artinian module,
or, equivalently, just consider the class of that module in the Grothendieck
group of the cartesian product of the localizations modulo any sufficiently high power of the radical.
\m
We cannot assemble this information into a homomorphism of Grothendieck
groups, because extracting the maximal Artinian submodule of a module
is not an exact functor.
This is how primary decomposition departs from the theory of algebraic cycles even on smooth quasiprojective varieties.
\m
In that special case, we can resolve the structure sheaf of a closed
subscheme by a finite sequence of locally free coherent sheaves.
Although  the cohomology of the constant sheaf ${\cal Z}$ could have torsion
one can use characteristic classes to recover the equivalence
class of a subscheme up to torsion and linear equivalence from the isomorphism
types of the terms of the resolution. Let's discuss this.
\vfill\eject\noindent 
{\bf Algebraic Cycles}
\m
The  cohomology and algebraic cycles of smooth quasi-projective
varieties are closely related to each other. Just as, in
algebraic topology, becuase of the Thom collapse, we know that two smooth real manifolds are cobordant
if and only if the Stieffel-Whitney numbers agree, it is also true
that if two divisors on a smooth complex projective variety are
linearly equivalent, they must have the same Chern numbers.
\m
In the appendix to the Borel-Serre paper on Riemann-Roch,
Grothendieck deduced a converse of this in every dimension, that two algebraic $i$ cycles
for any value of $i$ (not only $i=n-1$) are linearly equivalent modulo
torsion
on an $n$ dimensional smooth quasi-projective manifold if and only
if all the Chern classes (of the  of the
structure sheaves viewed as coherent sheaves on the smooth manifold) are the same, and that it is equivalent
to the difference between the elements in the Grothendieck group
being a torsion element plus a cycle of lower dimension.
\m
There, the structure sheaf of a sub-variety is resolved by
section sheaves of vector-bundles, and the Chern classes of a positive
vector bundle are defined to be the scheme of linear dependence of
each number of global sections.  
\vfill\eject\noindent
{\bf Lipman's question/conjecture about Chern cycles}
\m
The cohomology types of the Chern  subvarieties come from Bruhat
cells. While their cohomology classes,  the Chern classes
of vector bundles,  are  understood that way algebraically,
there remain are deep mysteries about what they actually look like geometrically.
Even if
 $M$ is  complex affine space  of dimension $n$ and  $V$ is its tangent
bundle, the scheme where vector-fields $v_0,...,v_s$
fail to be linearly independent,
has irreducible components all of dimension at least $s.$ The
deep and unsolved question of Lipman, formulated when he was a student
of Oscar Zariski, determines
 whether
there is actually a component of dimension  at least $s$  in every algebraic and  non-smooth  leaf closure, within the locus
where the $v_i$ are involutive, meaning $[v_i,v_j] =\sum a^{i,j}_k v_k$
And
it seems likely that algebraicity is not relevant for the question.
\m
The  smooth locus is where the
subsheaf of ${\cal O}_M^{s+1}$ spanned by the
action of directional derivative in each coordinate $(v_0(f),...,v_s(f))$
for $f\in {\cal O}_M$ local sections,
is locally free in its own right. This just describes the points where  we can enlarge
the span of $v_0,...,v_s$ by a local meromorphic change of basis to make
them independent anyway.
\vfill\eject\noindent
{\bf The total Chern class and cohomological definition}
\m
The so-called total Chern
class is the sum of the Chern classes, and it transforms direct sums
of vector bundles into cup products of cohomology classes.
\m
For a finite set of points of the Riemann sphere,
we choose a rational function $f$ sending these to a point $\infty$
of the Riemann sphere, and there is a complementary set of points
sent to another point $0.$ If we separate the two sets by a contour
$C,$ the  integral $\int_C {{df}\over f}$ can be interpreted as
$2\pi i$ times the number of poles on one side, or zeros on the other.
\m
The one-forms which have poles matching ${{df}\over f}$ and no zeroes
form a torsor, described by the additive Cech cocycle which is ${{df}\over f}$
restricted to a neighbourhood of the contour $C,$
and the same integral is an isomorphism invariant for this torsor.
If we call the finite set with its subscheme structure  $S,$ it is also the transform
along ${\cal O}_{{\mathbb P}^1}\to {\cal O}_S$ of the extension class of the Poincare
residue sequence $0\to \Omega_{{\mathbb P}^1} \to \Omega_{{\mathbb P}^1}(log\ C)\to {\cal O}_S\to 0$ in $Ext^1({\cal O}_S,\Omega_{{\mathbb P}^1}).$
\vfill\eject\noindent
{\bf Chern character}
\m
To define a  Chern character of an algebraic cycle  in this setting
it suffices  to define the Chern character of a locally free sheaf,
and by ignoring torsion  nothing else is lost tensoring cohomology with ${\mathbb Q}$
and defining Chern classes in the usual way
as vanishing of tuples of sections (after twisting by a line bundle)
and defining the Chern character as a rational polynomial function of
the Chern classes. The resulting cohomology class in degree $i$ for a locally free coherent
sheaf ${\cal F}$ if we bypass Chern classes altogether, doesn't require any denominators in its description, nor does it require reducing modulo
torsion; it is merely the one represented by the Hochschild $i$- cocycle
of $({\cal F},{\cal F}\otimes \Omega_M)$
sending local secions $x_1,...,x_i$ to the function sending
local sections $f$ of ${\cal F}$ to $f\otimes dx_1\wedge...\wedge dx_i.$
\m
The passage from first Chern classes to Chern characters is by
splicing resolutions, and there have been two attempts I know of
to generalize to the singular case. One thing a person could do is
consider non-exact sequences of locally free sheaves (the theory
of derived categories).
\m Here, such `chain complexes' are locally  objects that correspond
to homotopy types of semisimplicial modules (by the Dold correspondence [21]).
However it isn't clear why that homotopy type
ought to reflect meaningful properties of an algebraic cycle.

\m
The same cohomology class (of type $(1,1)$  is the extension class
of the exact sequence $0\to L\otimes \Omega \to {\cal P}(L) \to L\to 0$
twisted by $L^{-1}$ and it is induced from the extension
class
of the  Poincare residue
sequence itself 
$0\to \Omega \to \Omega(log D)\to {\cal O}_D \to  0$ by the
surjection ${\cal O}_M\to {\cal O}_D.$
The class has a homological definition as does the Chern character,
which therefore can be defined directly and `integrally' without
needing recourse to the historically earlier notion of Chern classes
of vector bundles.  
Even  integrally, Grothendieck's appendix
applies unchanged, although torsion in the Grothendieck group is not detected
by Chern character with coefficients in ${\mathbb Q}, $ nor is
it correctly detected if we use constant coefficients ${\mathbb Z}$ or
any other constant coefficients.
\vfill\eject\noindent
{\bf Homotopy theory}
\m
One might question whether the techniques of 
resolutions and cohomology are best; in the theory of derived
categories, one  interprets a resolution by locally free sheaves
as not only a formal alternating sum in the Grothendieck group, but
one remembers the chain homotopy type of the resolution. This is more
information as I mentioned; very precisely, on each affine part it captures the
homotpy type of a corresponding semisimplicial module by the
Dold correspondence. However it isn't clear why that homotopy type
ought to reflect meaningful properties of an algebraic cycle.
\m
The notion of considering the structure sheaf of a closed subscheme
as a coherent sheaf on an ambient variety and resolving it 
applies when things are nonsingular; if the ambient variety is singular
even when one can still find the necessary resolutions, or when
one allows non-exact resolutions, 
the particular step where a Koszul  resolution of the conormal sheaf
contributes the denominator of a Todd class would need to be modified
since the free resolution would not be a Koszul resolution.
\m
An attempt to reconcile Chern classes with the
singular setting may require considering that long exact sequences
resemble filtrations; in homotopy theory one can to some extent re-assemble
long exact sequences into exact couples, and this was useful for
determining homotopy types, especially stable homotopy types. 
In the theory of coherent sheaves, just passing to long exact
sequences captures for filtrations only successive pairwise extensions
and a little more than that, but misses the re-assembling that takes
place in a more correct analogy with homotopy theory.
The difficulty is caused by the fact that we define
coherent sheaves by keeping track of sections over open sets rather than
constructible sets. The denomninator in the Todd genus comes from
an assumption of local complete intersection, for example, and to get
past that one ought to be able to iterate Poincare residues when
nesting is more general; this requires considering sections on
constructible sets which are not necessarily open. We will not deal with that here.
\vfill\eject\noindent
{\bf Valuations}
\m
In conclusion, extending the theory of algebraic cycles  past the smooth case requires abandoning
the Todd genus in Riemann-Roch, whose denominator is the
class of the exterior algebra of the conormal sheaf of an algebraic cycle.
It requires abandoning the particular formula for the Todd genus, 
since the exterior algebra of a conormal sheaf need not be resolving in the non-smooth case. It requires considering that chain complexes of coherent
sheaves only capture something like the abelianization of a homotopy type,
and finally that coherent sheaves themselves as objects of a category
have a serious limitation when one looks at examples of trying to 
recapture algebraic cycles from their Chern classes, that one would need
to consider sections over constructible sets rather than only over open sets.
\m
A theory of valuations is a first attempt to rescue what one can
in the singular case.
\m {\bf Integral closure}
\m
If we make no distinction between saying $f$ is a local section of ${\cal F}$
versus $f{\cal G}$ agrees with a subsheaf of ${\cal F}{\cal G}$, the
tensor product modulo torsion, for ${\cal G}$ torsion-free of rank one,
on an irreducible scheme, then we are saying that we make no distinction
between local sections of ${\cal F}$ versus local sections  of its `integral closure.'
\m
If we {\it need not} make any such distinction, for ${\cal F}$ the structure sheaf, we are saying rational maps to ${\mathbb P}^1$ have no
a meromorphic section of a line bundle defined in codimension two
extends to a legitimate meromorphic section.
\vfill\eject\noindent
{\bf Review of construction of ${\mathbb R}$}
\m
One way to construct the reals ${\mathbb R}$ is
as contravariant functors ${\mathbb Q}\to Sets,$
an example of the now old construction which one can find in [22], and works with the ${\mathbb Q}$ as a partially-ordered set replaced
by any small category. According to Yoneda's lemma We identify a rational number $\alpha$
with the representable functor, the closed interval $[-,\alpha]$
leaving the first argument unspecified. Then for example
the colimit of $[- , \alpha]$ over $\alpha$ such that $\alpha^2<2$
is the functor $[-,\sqrt{2})$ which assigns to any rational number
$\beta$ the rational  points in the interval $(\beta,\sqrt{2}).$
\m
{\bf A group completion}
\m
For a reduced and irreducible scheme $M$ let's make a definition
which surely must pre-date Weil divisors; it agrees with the Weil
divisor group for smooth projective curves, but in general is larger.
\m
We begin with the standard notion  that  a {\it fractional ideal sheaf}
is
 a torsion-free coherent integrally closed rank one sheaf ${\cal I}$
labelled by a nontrivial  map (embedding)  ${\cal I}\to K$  of quasi-coherent sheaves to the constant
sheaf of rational functions
 $K.$
The fractional ideal sheaves  have a binary operation of multiplication
which sends ${\cal I}$ and ${\cal J}$ to ${\cal I}\otimes {\cal J}/torsion$ labelled
with the product labelling. It is the natural image of ${\cal I}\otimes {\cal J}$ under the tensor product of the two embeddings. However,
for our purposes we need a different binary operation on only the
set of integrally closed fractional ideals.
\m
We say  a fractional ideal sheaf ${\cal I}$  is {\it integrally closed} if
there is no strictly larger fractional ideal  sheaf ${\cal I}'$ with fractional ideal ${\cal J}$
such that ${\cal J}{\cal I}={\cal J}{\cal I}'.$   
\m Now we introduce a binary operation that is different than just
product of fractional ideal sheaves
\m
{\bf Lemma.} The binary operation
on integrally-closed fractional ideal sheaves assigning
to ${\cal I}$ and ${\cal J}$ the integral closure of
${\cal I} {\cal J}$ is associative and satisfies
cancellation.
\m
{Proof.} The only thing needing proof is the cancellation rule.
Suppose that  ${\cal I}{\cal J}$ and ${\cal I}{\cal L}$ have the same
integral closure where ${\cal J}$ and ${\cal L}$ are integrally
closed. Then there is some ${\cal M}$ so that ${\cal M}{\cal I}{\cal J}={\cal M}{\cal I}{\cal L}$ and so taking the product with ${\cal M}{\cal I}$ makes
${\cal L}$ and ${\cal J}$ become equal, by definition of integral closure
they are already equal.
\m
It follows that
\m
{\bf Corollary.} The monoid of integrally closed fractional ideal sheaves
with this binary operation embeds faithfully in to its group completion $\Gamma.$.
\m
We denote by $\nu({\cal I})$ the image of the integrally closed
fractional ideal ${\cal I}$ in this group. Although the letter $\nu$
corresponds linguistically to the roman $n,$ its visual similarity
for $v$ means it has been used traditionally to stand for `valuation.'
\m
Also, if ${\cal I}$ is not integrally closed, we also write
$\nu({\cal I})$ for the image in $\Gamma$ of the integral closure of ${\cal I}$.
Thus we think of $\nu$ as a sort-of universal valuation, taking
values in the commutative group $\Gamma.$
\vfill\eject\noindent
{\bf `Rational' fractional ideals}
\m
Next, instead of choosing a character $\nu:\Gamma\to {\mathbb Q}$ 
which might describe a particular valuation,  we
should merely consider the natural $\Gamma \cong \Gamma\otimes {\mathbb Z}
\to  \Gamma\otimes {\mathbb Q}$ induced by ${\mathbb Z}\subset {\mathbb Q}$.
By a slight abuse of notation we'll also denote by $\nu({\cal I})$ the
image of (the integral closure of) a fractional ideal ${\cal I}$ in 
$\Gamma\otimes {\mathbb Q}, $ however, when we do that, 
it is no longer true that $\nu$ faithfully represents integrally closed
fractional ideals, as the map $\Gamma\to \Gamma\otimes{\mathbb Q}$ will
have nontrivial kernel when $\Gamma$ has torsion.
\m
Let's 
say that an element $\gamma\in \Gamma\otimes {\mathbb Q}$ belongs to the {\it effective cone}
if there is a nonzero natural number $n$ such that $n\gamma$ is
the image of an element in the effective cone of $\Gamma.$
\m There is a partial ordering on $\Gamma\otimes{\mathbb Q}$ such that
 $\nu({\cal I})\le \nu({\cal J})$ if and only if after 
modifying ${\cal I}$ or  ${\cal J}$  by addition of a torsion element,
there are representatives
of ${\cal I'}$ and ${\cal J'}$ in $\Gamma$ and a number $N$ so that
(the integral closure of) ${\cal J'}^N$ is contained in the integral closure
of ${\cal I'}^N.$
(Note that the inequality is reversed from
the familiar one that relates inclusions of sets with their cardinalities).
\m
Thus we think of $\nu$ as a sort-of universal valuation, taking
values in the commutative group $\Gamma\otimes {\mathbb Q}.$
\vfill\eject\noindent
\m
{\bf Completion of the group-completion}
\m Finally, we are going to complete $\Gamma\otimes {\mathbb Q}$ in the other sense.
The completion
of $\Gamma\otimes {\mathbb Q}$ can be defined using
representable functors, as we did for ${\mathbb Q}$ itself already,
I expect that the result will be no different than
$\Gamma\otimes {\mathbb R}.$
In that case
we are now working in a real vector-space $\Gamma\otimes{\mathbb R}$ with an effective
cone.

\vfill\eject\noindent
{\bf Nash transforms}
\m
In the case of iterating Nash transforms,  one has a sequence
of varieties $...V_2\to V_1 \to V_0=V$ and if we call $\Omega_{ij}$
the pullback modulo torsion to $V_i$ of the one-forms of $V_j,$ 
then taking $n$'th exterior power and again reducing modulo torsion
we obtain a system of locally free sheaves, which we can write
${\cal O}_{V_i}(K_j)$ with $K_j$ a Cartier divisor on $V_j$ for $i>j.$
\m
One thing which we will do later, 
as naturally as possible, is to  restrict to a subvariety of $V$
which supports $n$ generically independent vector-fields, and to describe ${\cal O}_{V_i}(K_j)$ 
for $j<i$ as an ascending series of  locally principal  ideal sheaves 
on $V_i$ using the theory of connections. On the union of the regular locus $U_0\to U_1\to U_2...$
there is then an infinite ascending chain of ideal sheaves, 
and the colimit is a locally principal ideal sheaf on the whole
of the Nash manifold $U.$ 
\m
The ideal sheaves are not `relatively basepoint-free' however
when $V_0$ and therefore all $V_i$ are irreducible, 
any  fractional ideal sheaf on any $V_i$ 
 can be described as the tensor difference
of two ideal sheaves on $V_0.$ A fractional ideal sheaf is just
a torsion-free rank-one sheaf labelled by an embedding in the 
constant sheaf of rational functions, and the same is true
more simply if we ignore labellings, that any coherent sheaf on any $V_i$  
is the result of pulling back from $V_0$ a torsion free coherent sheaf
whose reduction modulo torsion becomes invertible, and tensoring the
inverse with the pullback of another coherent sheaf on $V_0$, tensoring
and again reducing modulo torsion.
\m
In this sense, the $i$'th locally principal ideal sheaf on $U$ in our
series  is the result of pulling back $\Lambda^{n+1}{\cal P}({\cal F}),$
tensoring with the inverse of the pullback modulo torsion of ${\cal F}^{r+1}$
and again reducing modulo torsion.
And here ${\cal F}$ must be defined recursively, to be the 
product of the first, second, third, ... $(i-1)$'st  locally principal
ideal sheaves. 
\vfill\eject\noindent
Because of naturality of first principal parts, we may apply ${\cal P}$
at any stage on or past the $i$'th stage. In some sense the
most efficient way is to apply it always at the $i$'th stage, and then
the result is not yet invertible, all the earlier ones are, and
in fact the $i+1$'st Nash transform is the universal one which makes
this invertible.
\vfill\eject\noindent 
{\bf A natural transformation}
\m
It is convenient to leave the torsion-free rank one sheaf ${\cal F}$
unspecified, and then there is a natural map (natural transformation 
of functors)
$$\Lambda^{n+1}{\cal P}({\cal F})^{\otimes n+2} \to
\Lambda^{n+1}{\cal P}({\cal F}\otimes \Lambda^{n+1}({\cal P}({\cal F})) \ \ (1)$$
which one can reduce modulo torsion. 
It is analgous to `carrying' when adding base $n+2$
expansions.
If we write 
${\cal F}_0={\cal O}_M,\ 
{\cal F}_1=\Lambda^{n+1}{\cal P}({\cal F})/torsion, \ {\cal F}_2={\cal F}_1^{\otimes 2}/torsion , \ ... \ {\cal F}_{n+1}={\cal F}_1^{\otimes (n+1)}/torsion,
{\cal F}_{n+2}=\Lambda^{n+1}{\cal P}({\cal F}{\cal F}_1)/torsion, {\cal F}_{n+3}={\cal F}_{n+2}{\cal F}_1/torsion,...$ and in general
${\cal F}_{(n+2)^{i+1}}=\Lambda^{n+1}{\cal P}({\cal F}{\cal F}_1{\cal F}_{(n+2)}...{\cal F}_{(n+2)^i}),$
then  from  the maps ${\cal F}_i^{\otimes n+2} \to {\cal F}_{i(n+2)}$ 
we can derive existence of particular natural   maps ${\cal F}_a\otimes {\cal F}_b \to {\cal F}_{a+b}/torsion.$
The ${\cal F}_i$ are torsion-free and rank one, and are generated
by global sections if ${\cal F}$ is very ample. 
\m
In [23] in an initial result, I showed that in characteistic zero the necessary and sufficient condition for
Nash's process to finish in finitely many steps is that starting
with ${\cal F}={\cal O}_M,$ 
eventually 
${\cal F}_{(n+2)^i}^{\otimes n+2}/torsion 
\to 
{\cal F}_{(n+2)^{i+1}}/torsion$ should become an isomorphism after multiplying by a power
of ${\cal F}_i$ and ${\cal F}_{i+1}.$
\m
The  argument  in [24] using finiteness of normalization improves this, showing that in fact
that the following simpler condition is necessary and sufficient, that
eventually. ${\cal F}_{(n+2)^i}^{\otimes n+2}\to {\cal F}_{(n+2)^{i+1}}$
is surjective (equivalently, an isomorphism).
\m
Therefore, a necessary and sufficient condition 
for the existence of a torsion-free coherent rank one resolving
sheaf ${\cal F}$ (whose existence is equivalent to Hironaka's theorem) is that on the first stage
$$\Lambda^n{\cal P}({\cal F})^{\otimes  (n+2)}/torsion \to 
\Lambda^{n+1}{\cal P}({\cal F}\otimes \Lambda^{n+1}({\cal P}({\cal F}))
$$
\m
Here is how the proof in [24] works. The weaker statement is equivalent where
we allow ourselves to  multiply by a power of ${\cal F}\Lambda^{n+1}{\cal P}({\cal F})/torsion$. The way to prove that multiplying by this power is
unnecessary 
is to assemble ${\cal F}_0\oplus {\cal F}_1 \oplus {\cal F}_2...$ into a sheaf of  graded rings on $V_0$ . The indexing of the ${\cal F}_i$
 reflects this grading which wasn't evident 
in earlier notation.
\m
Here is how we prove that  surjectivity 
after multiplying by a  power of ${\cal F}\Lambda^{n+1}{\cal P}({\cal F})$
implies we can find a new choice of ${\cal F}$ where surjectivity
holds without needing to multiply by a power like that.
The combinatorics of adding base $n+2$ expressions for natural numbers
imply that finite type of the graded sheaf of rings is equivalent to one of  the stronger equalities
holding, with ${\cal F}$ replaced by some ${\cal F}_1{\cal F}_{(n+2)}...{\cal F}_{(n+2)^i}={\cal F}_{{{(n+2)^{i+1}-1}\over{n+1} } },$ whereas the weaker condition
where we allow multiplying by the power only   implies that the
normalization of the graded sheaf of rings is finite type.
\m
We deduce the stronger equality from the weaker one, for a possibly
larger value of $i$, as follows: 
Let's describe the affine case so the language of sheaves doesn't get
in the way. The subring of the graded ring generated up to each
degree $i$ is finite type and has a normalization which is finitely
generated as a module over this ring. The series of normalizations
is an increasing series of subrings of the overall normalization,
which exhausts the overall normalization, then contains all of
the finite number of ring generators of that. It follows
that the normalization of one of the finite type subrings is the
same as the overall normalization. But the intermediate subrings
between these are an ascending chain of submodules of a finitely-generated
module over a Notherian ring. This completes the proof.
\vfill\eject\noindent 
{\bf Discussion}
\m
While there are no natural inclusions among the ${\cal F}_i,$
 there are
 natural maps 
${\cal F}_{(n+2)^i}^{\otimes n+2}/torsion 
\to {\cal F}_{(n+2)^{i+1}}.$
When ${\cal F}={\cal O}_M$ then ${\cal F}_i$ is naturally
a subsheaf of the rational function field tensor $\Lambda^n\Omega_M^{\otimes i}.$
If one wishes a premonition of our  arguments using connections, one can
merely calculate the ${\cal F}_i$ when we have embedded $\Omega_M$
into a free module using a sequence of generically independent
global derivations, and see that inside the fraction field
tensor $\Lambda^n\Omega_M^{\otimes i}$  lives inside 
of ${\cal O}_M.$

\m
Coherent sheaves make sense analytically, and finding an ${\cal F}$
which makes (1) surjective comprises a second order differential equation
and we've just shown that solving it  is equivalent to resolving the singularitis of $M.$ The equivalence,
ever, being indirect. 
\m
Let's re-describe what we've just done in a more geometric language.
Once we resolve the indeterminacies of ${\cal F}$
to obtain a $W_0\to M_0,$ where the pullback modulo torsion of ${\cal F}$
is invertible,  and then if we further resolve
indeterminacies of $n$-forms on $W_0$ to obtain $W_1\to W_0$
the pullback of ${\cal F}\otimes \Lambda^{n+1}{\cal P}({\cal F})$
modulo torsion becomes an invertible sheaf, 
is the one which  we were calling ${\cal F}{\cal F}_1$ above, 
with the mulitiplication being tensor product modulo torsion,
and our map is
$${\cal F}_1^{\otimes n+2} \to \Lambda^{n+1}{\cal P}({\cal F}{\cal F}_1).$$
Tensoring with the $-n-1$ tensor power of ${\cal F}_1$
then gives
$${\cal F}_1\to {\cal F}_1^{\otimes -n-1}\Lambda^{n+1}{\cal P}({\cal F}{\cal F}_1)$$
and ${\cal F}_1$ itself has a similar form to the right side of the
display; if we write the pullback modulo torsion of ${\cal F}$ itself
as ${\cal I}$ it is  ${\cal I}^{\otimes -n-1} \Lambda^{n+1}{\cal P}({\cal I}).$
So we are looking at a map
$${\cal I}^{\otimes -n-1}
\Lambda^{n+1}{\cal P}({\cal I})
\to 
({\cal F}{\cal F}_1)^{\otimes -n-1}\Lambda^{n+1}{\cal P}({\cal F}{\cal F}_1).$$
Both the domain and codomain are expressed as a multiplicative difference of relatively
basepoint-free coherent sheaves, and this map, on $W_1,$ is precisely
the one which determines ramification; it is the map from the
pullback modulo torsion of the $n$-forms of $W_1$ to the $n$-forms
modulo torsion of $W_2.$
\m
If the {\it original} map (1) is surjective, this one is too, and 
$W_1\to W_0$ is an unramified Nash transformation. Note that this
is not quite the same as saying that $W_1$ resolves the singularities
of $M_0.$  However, $W_1$ is nonsingular, because we are in characteristic zero,
or working analytically, and the map now expresses an isomorphism between
the $n$-forms modulo torsion and a locally free sheaf.
\m
Conversely, if ${\cal F}$ is a resolving sheaf, $W_0$ is smooth, 
and $(1)$ becomes surjective after pulling back and reducing modulo
torsion. This only implies that $(2)$ becomes surjective after
tensoring with a tensor power of ${\cal F}\otimes \Lambda^{n+1}{\cal P}({\cal F}),$ however, as we showed above,  that if we repeatedly replace
${\cal F}$ by ${\cal F}\otimes \Lambda^{n+1}{\cal P}({\cal F}) $
and the {\it just once} apply $\Lambda^{n+1}{\cal P}( -- )$
then when we replace ${\cal F}$ by the new sheaf, the map
(1) has become surjective. l
here since it uses Hilbert's basis theorem.
\m
It is interesting that all we have been talking about so far is
whether a second order differential equation can characterise 
resolutions, and it has led us into needing to consider the same
algebraic structure which relates to Nash's question.
\m
The discussion generalizes if we wanted to talk about either
algebraic cycles and cohomology, or indeterminacy of coherent rank one
sheaves, however let's focus here on the resolution problem.
\vfill\eject\noindent
{\bf Moduli of resolutions}
\m
The isomorphism types of torsion-free coherent ${\cal F}$ which make (1) surjective comprise a 
type of infinite-dimensional moduli
space.  Although the relation with resolving is  indirect
(on the one hand, requiring us to replace ${\cal F}$ by some higher ${\cal F}_i$if we wish to literally  
obtain a resolving sheaf, and conversely requiring us to replace
the purported resolution by a Nash blowup we wish it to literally  correspond toone of the  ${\cal F}$), a way of approaching
the question of finding resolutions within a finite-dimensional moduli space would be to filter 
moduli of solutions of the surjectivity of (1) by finite-dimensional spaces.
\vfill\eject\noindent
{\bf Connections}
\m Let's suppose we're given a
connection 
$$\nabla:{\cal L}\to {\cal L}\otimes \Omega$$
satisfying $$\nabla(fy)=f\nabla(y) +y \otimes df$$
for local sections $y$ of ${\cal L}$ and $f$ of ${\cal O}_M.$
\m
The fibre vectoriel $L$ whose structure sheaf pushes down on $M$
to ${\cal O}_M\oplus {\cal L} \oplus ....$ 
has one-forms 
$$({\cal O}_Md{\cal O}_M) \oplus ({\cal L}d{\cal O}_M + {\cal O}_M d'{\cal L}) \oplus ...$$
Here $d'$ is the differential on the line bundle, and when
we take  forms with logarithmic ppoles on the zero section
and twist  by the zero section given multiplicity $-1,$  The `degree-one' term which survives is ${\cal P}({\cal L}).$
The first summand is a subsheaf isomorphic with  ${\cal L}\otimes \Omega_M \oplus {\cal L},$
while the second summand is only a coherent subsheaf after we reduce
modulo the first term, since $d$ is not ${\cal O}_M$ linear.
However, the set of local sections $-\nabla(y)+1\otimes dy$ does describe
a subsheaf, as  for $y$ a local section of ${\cal L}$ we have $-\nabla(fy)+d(fy)=-f\nabla(y)-y \otimes df + y\otimes df 
+f\otimes dy  = f\cdot (-\nabla(y)+dy)$ for local sections $f$ of
${\cal O}_M $ showing it is stable
under multiplication by  local sections of ${\cal O}_M.$
\m
Locally, once we choose a generating section $x$ of ${\cal L}$
then for each $f$ we have
$$\nabla(yf)=y\nabla(f)\oplus f\otimes dy.$$
\m
In fact, ${\cal L}^{-1}{\cal F}$ is an ideal sheaf,
and $y$ ranges over local sections of this ideal sheaf, and 
we have the element $\nabla(f)$ depending on our choice of generator,
and locally ${\cal P}({\cal F})$ is isomorphic with the subsheaf of
${\cal F}\oplus \Omega$ spanned by $y\oplus dy.$ The
isomorphism  depends on the choice of $f,$ and it sends
$y\oplus dy$ to $y\nabla(f)\oplus f\otimes dy.$ 
\m
Therefore we could define ${\cal P}({\cal F})$ by a type
of additive cocycle.  If we call the ideal sheaf ${\cal I}$
first we construct on each open set where $f$ is a local
section the subsheaf of ${\cal I}\oplus \Omega$ spanned
by the $y\oplus dy.$
\m
Covering $M$ by open sets, 
the patching Cech cocycle for  subsheaf described in the summand on the
right side is easy to describe, it is just the ordinary patching in ${\cal F}\otimes \Omega.$ On
an open set where $f$ is a chosen generator of ${\cal L}$ and another
where $g$ is a chosen generator, the patching isomorphism is
multiplication by the unit  $u=g^{-1}f.$
\m
For the patching of the whole sheaf,  when $g=uf$ as before,
each $\nabla(f)$ is equal to 
$$\nabla(ug)=u\nabla(g)\oplus g\otimes du
$$ $$ =f \ g^{-1}\nabla(g) + g\otimes (gdf-fdg)/g^2$$
$$=f g^{-1}\nabla(g) 
+ df- f {{dg}\over g}$$
\m
The nicest way of writing this is in terms
of   $\nabla \ log(g) = {1\over f}\lambda (f),$ where it says
$$\nabla log(f)-\nabla log(g) = dlog(f)-dlog(g).$$
\m
Instead of saying that  ${\cal P}({\cal F})$
is spanned by the $y\nabla(f)\oplus f\otimes dy$
we could have said that it  is spanned by the $yf (\nabla(log\ f) \oplus  1 \otimes d\ log(y)).$
\m
If we multiply $f$ by a local unit and divide $y$ by the same,
if we call this unit $s,$ then $nabla(log \ f)$ becomes
${1\over {sf}}
\nabla(sf)=
{1\over{sf}}
(s\nabla (f) + f\otimes ds)
=\nabla\ log(f) + d log (s) $
while $d\ log(y)$ becomes $d\ log(y)-d\ log(s)$  and the expression
is unaffected.
\m
Therefore the element $(\nabla\ log \ f \oplus 1\otimes d\ log \ y)$
is a well-defined rational section of first principal parts depending only
on the product $yf,$ and if we had named the spanning element
with the symbol $\nabla'(yf)$ with $\nabla'$ denoting what we would
think of as a  universal connection, we would be factorizing $\nabla'(yf)=yf \nabla'\ log \ yf.$
\m
The universal connection $\nabla'$ can likely be understood directly by thinking of first principal
parts related to the square of the ideal defining the diagonal.
\m
In an earlier explicit description we said when we have
generically independent  derivations $\delta_1,..,\delta_n$ with 
values in the rational function field $K$
then we define ${\cal P}({\cal I})$ for ${\cal I }$ an ideal
sheaf, or fractional ideal sheaf,  to be spanned by local sections  $y\oplus \delta_1(y) \oplus...
\oplus \delta_n(y)$ in ${\cal F}\oplus  K^n$ for $y$ local
sections of the ideal sheaf.
\m
Recall that $(\nabla\ log \ f \oplus 1\otimes d\ log \ y)$
depends only on $yf$ for $f$ a local section of ${\cal L},$
and when $f=1$ this is just $1\otimes d \ log \ y$
So the product with $y$ becomes $y\oplus dy$, and 
we see that the description of the universal connection
agrees with what we had said when we described principal parts
using derivations.
\m
Earlier I wrote that if ${\cal I}$ is spanned
locally by secions $t_i$ and ${\cal O}$ is spanned
by $x_i$ then letting $s_i$ be the products
$t_jx_k$   $\Lambda^{n+1}{\cal P}({\cal I})/torsion$
is spanned by $(s_0\oplus \nabla s_0)\wedge...\wedge (s_n\oplus\nabla(s_n))$
where $1\oplus \nabla$ is the `universal connection,'
and this calculated to $s_0...s_n dlog (s_1/s_0)\wedge...\wedge dlog(s_n/s_0).$
\m
By $1\oplus \nabla$ I really meant that  we can calculate the universal
and invaraint connection 
$\nabla'$ in terms of our ordinary connection $\nabla$ 
using our more rigorous formula 
$\nabla' \ log \ yf = \nabla \ log \ f + d \ log \ y$
with $\nabla$ our actual connection
\m
Thus there is a theoretical way of describing the original
construction using derivations as coming from a choice of
connection, but in any case when ${\cal F}$ is an ideal
sheaf and $\delta_i$ derivations of ${\cal O}_M$ 
generically independent and well-defined on a neighbourhood
we can use that to embed ${\cal P}({\cal F})$ into ${\cal O}_M^{n+1}$
.
\m
For some reason we wanted to take the $\delta_i$ commuting.
\m
We always have
$${\cal F}_{(n+2)^\alpha}^{n+2}\subset {\cal F}_{(n+2)^{\alpha+1}}
\subset ({\cal F}{\cal F}_1...{\cal F}_{(n+2)^\alpha})^{n+1}$$
and if I replace the first factor ${\cal F}^{n+1}$ on the right
with ${\cal F}_1$ we can compose the chain of inclusions
${\cal F}_1^{n+2}({\cal F}_{(n+2)}...{\cal F}_{(n+2)^\alpha})^{r+1}\subset {\cal F}_{(n+2)}^{n+2}({\cal F}_{(n+2)^2}...{\cal F}_{(n+2)^\alpha})^{r+1}\subset ... \subset {\cal F}_{(n+2)^\alpha}^{r+2}.$
\m
Now recall that for a fractional ideal ${\cal I}$ we denote by $\nu(I)$ the associated element of $\Gamma\otimes{\mathbb R}$
\m
We are going to apply this to the chain of inclusions we've established
$$...{1\over{(n+2)^\alpha}}\nu({\cal F}_{(n+2)^\alpha})
\ge
{1\over 
{(n+2)^{\alpha+1}}
}
\nu({\cal F}_{(n+1)^{\alpha+1}})
\ge ...  $$ $$\ \ \ \ \ \ \ \ \ \ ... \ge 
{{n+1}\over{(n+2)^{\alpha+1}}}\nu({\cal F}{\cal F}_1{\cal F}_{(n+2)}
...
{\cal F}_{(n+2)^\alpha}) 
 \ge{{n+1}\over{(n+2)^\alpha}} \nu({\cal F}{\cal F}_1...{\cal F}_{(n+2)^{\alpha+1}})
$$ $$\ge   {1\over{(n+2)^\alpha}}\nu({\cal F}_{(n+2)^\alpha})
-{{n+1}\over{(n+2)^\alpha}}\nu({\cal F})+{1\over{(n+2)^\alpha}}\nu({\cal F}_1).
$$
\m I'd like to explain the somewhat unusual use of elipses ... here.
Except for the last term written one sees in infinite increasing chain
followed by an infinite decreasing chain of larger numbers. The last inequality written
then shows that  each term of the decreasing chain of larger numbers can be 
made less than  corresponding element of the smaller increasing chain
by subtracting  a  rational multiple of 
the constant $\nu({\cal F}_1)-(n+1)\nu({\cal F}),$ and the coefficients
in that rational multiple tend to zero.
\m
Each fractional ideal is a functor of ${\cal F}$ of a particular
degree, and the denominator we've introduced as a coefficient to
each $\nu$ value is the degree.  That is to say, if we abbreviate
by $\nu({\cal F})$ the expression ${1\over{degree(F)}}\nu({\cal F})$
then our inequalities become 
$$...\mu({\cal F}_{(n+2)^\alpha})
\ge
 \mu({\cal F}_{(n+1)^{\alpha+1}})
\ge ... $$
$$\ \ \ \ \ \ ... \ge
(n+1)\mu({\cal F}{\cal F}_1{\cal F}_{(n+2)}
...{\cal F}_{(n+2)^\alpha})
\ge (n+1) \mu({\cal F}{\cal F}_1...{\cal F}_{(n+2)^{\alpha+1}})
$$ $$ \ge   \mu({\cal F}_{(n+2)^\alpha})
-{{n+1}\over{(n+2)^\alpha}}\nu({\cal F})+{1\over{(n+2)^\alpha}}\nu({\cal F}_1).
$$
It follows that the $\mu({\cal F}_{(n+2)^\alpha})$ form a decreasing
sequence, and the $(n+1) \mu({\cal F}{\cal F}_1...{\cal F}_{(n+2)^\alpha})$
comprise an increasing sequence, and 
\m
{\bf Theorem.} Provided a suitable connection exists (i.e. locally
on $M$)  then
the limit of the increasing sequence in the
completion of $\Gamma\otimes {\mathbb Q}$
is equal to  the limit of the decreasing sequence. 
\m
The  fact that the limit of the sequence of interest, the increasing sequence,
is equal to the limit of the decreasing sequence means we may
change focus and try to understand the decreasing sequence.
\m Although, the definitions of the two sequences are intertwined
in a particular recursive way, also.
\vfill\eject\noindent
{\bf Pre-conclusion}
\m
My thinking on these matters is not independent.
M. Thaddeus, in a letter to Teissier about the  preprint with Perez, said that a purported open cover of a toric
variety actually covers only the rational points. The thinking above is partly generated
by this idea. For a toric variety, a rational character of the lattice corresponds to 
an equivariant rational curve, and the question whether minimum point of the  image of our ascending sequence
stabilizes is equivalent to lifting the arc to the smooth locus of the Nash manifold. One
visualizes the arc lifted part-by-part extending across an infinite union of open subsets
exhausting the Nash manifold, and whether the open cover is finite.
\m
Although Thaddeus didn't specify what he means by a rational point, earlier,
speaking generally, M. Spivakvski said [10] that valuations
ought to be viewed as generalizations of arcs, and arcs can be
viewed as homomorphisms to a formal power series ring.
\m
My construction of $\Gamma$
and of a completion  of
$\Gamma\otimes {\mathbb Q}$  is meant to be able to speak
of a type of non-rational point which I think that both conversations
were referring to in different ways.  If we manage to skip the step
of tensoring with ${\mathbb Q}$ the completion might have
torsion. We will have a toric interlude,
and the last of this paper will establish the conditions
mentioned in [26] necessary and sufficient to lift a {\it holomorphic}
arc, although, our results above and in the  next section suggest
that lifting of arcs which go to ideal points of $M$ (in the sense that
the imaginary axis goes to a cusp in a modular curve),  would also be necessary to prove
properness of Nash's manifold.
\vfill\eject\noindent
{\bf Toric interlude}
\m
When $M$ is a projective toric variety and $\Lambda$ is the lattice of characters, 
once ${\cal F}$ is a very ample invertible sheaf generated by a finite set of
torus characters all ${\cal F}_i$ are also generated by a finite set of
torus characters. The generators are described in [24] Theorem 9 which uses [25] proposition 11.
The $\mu({\cal F}_i)$ can  be understood
as corresponding to particular convex subsets of either $\Lambda \otimes {\mathbb Z}[{1\over {n+2}}]$ or more simply, just $\Lambda \otimes {\mathbb Q}$. Each
$\mu({\cal F}_{(n+2)^{i+1}})$ corresponds precisely to the convex
hull of the set of averages of $n+1$ torus-equivariant global sections of
${\cal F}{\cal F}_{{(n+2)^{i+1}-1}\over{n+1}}={\cal F}{\cal F}_1...{\cal F}_{(n+2)^i}$ which are not on a hyperplane in $\Lambda \otimes {\mathbb Q}$
while each ${\cal F}{\cal F}_1 ... {\cal F}_{(n+2)^i}$ corresponds to 
the convex hull of the result of iterating averages of the larger number of $(n+2)$ 
elements not on a hyperplane, iterated $i$ times. Note in both cases
we do not pass to convex hull until after the iterated averaging is done. 
\m {\bf Remark.} 
The torus-equivariant global sections
of a coherent sheaf ${\cal G}$ generated in this way by a finite set  $T$  can be larger than $T,$
including also the elements which  conduct global sections of ${\cal F}^N$ into ${\cal F}^N{\cal G}$
for any $N.$ This is in general strictly intermediate between $T$ and its convex hull.
The inductive construction of [24] and [25] allows ignoring such extra generators 
to construct a correct generating set for each ${\cal F}_i$ without needing to 
know the full set of equivariant global sections in previous stages.

\m
A character $\Lambda \to {\mathbb R}$ can be evaluated on each
$\mu({\cal F}_i)\in \Gamma\otimes {\mathbb R}$ by evaluating on the corresponding subset of
$\Lambda\otimes {\mathbb Q}$ and choosing the minimum value. 
\m
At the same time, characters  $\Lambda \to {\mathbb R}$ 
correspond to arcs in the toric variety. Characters which
factorize through ${\mathbb Q}$  correspond to equivariant
rational curves, while
those which do not  correspond to a type of non-holomorphic arc. In both
cases
liftability to the smooth locus is determined by whether the increasing sequence
has finitely many distinct values in the image.
\m
If our basic torus characters are $x_1,...,x_n$ then a rational character determines in a primitive
generator a primitive rational monomial in the $x_i,$ while a real character can determine
a function of the  
type $x_i(t)=e^{\lambda_i t}$ with $\lambda_i$ positive real, which converge
not as $t\to 0$ but as $t \to -\infty$ along the real line.
\m
{\bf Example} When the equivariant global sections
of ${\cal F}$ correspond to  $\{(0,0),(0,1),(1,0),(1,3),(2,5)\} \subset {\mathbb Z}^2$ the figure on the next  page shows  convex hull of ${1\over {16}}$ times the
lattice point set $\{(5,5),(5,21),...\}$ corresponding to ${\cal F}_{16}$ and the convex hull
of ${1\over{48}}$ times the  lattice point set $\{(16,16),(16,64),...\}$ corresponding to 
${\cal F}_{48}.$ The theorem of the section above called `Connectoins' implies
that  as either of the powers $q,s$ is increased, the set-theoretic difference
of ${1\over {4^q}}$ times the hull of the lattice set of ${\cal F}_{4^q}$ minus
${1\over {3\cdot 4^s}}$ times the hull of the lattice set of ${\cal F}_{3\cdot 4^s}$ decreases monotonically or stays constant, and as  $q,s\to \infty$ the difference 
shrinks to a one-dimensional continuous curve of area zero.  For general projective
toric varieties the theorem implise that it will always be a topological sphere. This
is consistent with the stronger notion that the limiting curve is the boundary of a
polyhedron. For projective toric varieties, the stronger fact  will hold precisely when the inner curve
stabilizes in finitely many steps (the outer curve never does, but always converges to the limit of the inner curve). The stronger assertion remains unproven and would be
strictly stronger if and only if the arc lifting problem for such exponential types
of arcs, converging as $t\to -\infty,$  were inequivalent to the arc lifting problem for holomorphic
arcs, of the type converging when $t\to 0.$

\noindent 
\scalebox{.8}{\includegraphics {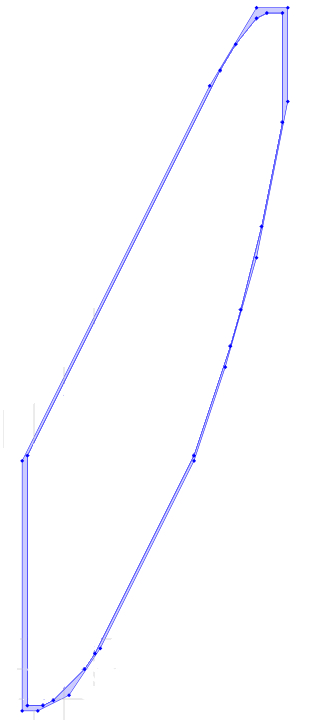} }

\noindent The javascript for graphing the convex hulls on a canvas is due to indy356.github.com.

\vfill\eject

\pagenumbering{arabic}
\noindent 
{\bf Proof of the theorem of [26]}.
\m
Let $M$ be a singular normal irreducible complex manifold of
dimension $n.$ 
Let $f:U\to M$ be the smooth Nash manifold, which is the smooth
locus of the inverse image of the tower of Nash transforms of $M.$
Let $\gamma:R\to M$ be an
analytic map from a  connected Riemann surface whose
image contains a nonsingular point of $M.$
\m
Conjecturally, it is sometimes said,
  there is a  unique $\gamma':R\to U$
such that $\gamma=f\circ \gamma'.$
\m
The Nash manifold $U$ supports  Cartier divisors $K_0,K_1,....$
so that the coherent sheaf ${\cal O}_U(K_i)$ is the pullback modulo torsion of
the differential $n$-forms of the $i$'th Nash transform; these are
allowed to have infinitely many irreducible components.
By a triangular change of basis, we let
$$D_{(n+2)^j}=K_j+(n+1)\sum_{t=0}^{j-1}(r+2)^{j-1-t}K_j$$ 
(we will define $D_j$ for non-powers of $n+2$ a bit later). The
reason for doing that is that
the divisors $D_j$ are relatively basepoint-free; in fact
we can give an elementary description of a coherent sheaf ${\cal F}_j$
on $M$ which pulls back modulo torsion to ${\cal O}_U(D_j).$
\m
One way of remembering the formulas above is  to solve for $K_U$
as an infinite series
$${{-1}\over {(n+1)}}K_U=D_1 + D_{(n+2)}+D_{(n+2)^2}+...$$
\m
The left side of the equation describes an effective Cartier divisor
now, albeit a ${\mathbb Q}$-divisor. If our arc is unbounded
the tails of the series would be allowed to get smaller
as we pass from consideration of one singular point to another,
and would be allowed to intersect trivially.
\m
The purpose of this note is to explain a theorem which is stated
in the earlier preprint, which in this language says
\m {\bf  Theorem.} Let $\gamma:R\to M$ be an analytic map
from a connected Riemann surface whose image contains just one non-singular
point of $M.$ Then there is a lift $\gamma':R\to U$
such that $\gamma=f \circ \gamma'$  if and only 
if the series of intersection numbers 
$$\gamma\cdot D_{n+2}+\gamma\cdot D_{(n+2)^2} + \gamma \cdot D_{(n+2)^3}...$$
is eventually a  (necessarily nonzero) geometric series to the base of $(n+2).$
\m
The theorem will extend to the case when the image of $\gamma$ contains
finitely many singular points; however, if $\gamma$ is unbounded
it is possible that the tails for the separate singular points might
intersect trivially; in that case we could apply the theorem to
each singular point to deduce arc lifting even while the
entire series would not need to be geometric.

\m
If we start with an arc $\gamma$ which is already lifted to $U$,
the case when the series converges literally is when it has only finitely
many nonzero terms. In this case, the arc maps to a single point of
$M$; it always will share its fiber with a complete curve
of the other type.  
\m
When the condition is holds with infinitely many nonzero terms,
and $\gamma$ meeting only one singular point,  it is even in this case not only an abstract condition.
The formula for the sum of a geometric series remains valid, and it is allowed to 
to be interpreted  formally or $p$-adically for a prime divisor of $n+2$.

\m The proof of the theorem as we've stated it here is
a formal manipulation.
A  bit later we will describe generators of
coherent ideal sheaves on $M$ itself which pull back modulo torsion
to the ${\cal O}_U(D_i).$ The proof is just formal, and  all we need to know is
that all are effective. 
The formal proof is to rewrite the right side as
{\normalsize
$${{1}\over{n+1}}( -D_1 + ((n+2)D_1-D_{(n+2)}) + ((n+2)D_{(n+2)}-D_{(n+2)^2})+....)$$
}Multipying both sides by $n+1$ gives
{\normalsize $$-K_U=-D_1 +((n+2)D_1-D_{(n+2)})+((n+2)D_{(n+2)}-D_{(n+2)^2})+.....$$
}
\m
Each term in round parentheses is one of the  effective
divisors $-(K_{i+1}-K_i)$ 
\m It does not matter whether we view these as divisors
on the Nash manifold itself or, term-by-term, each on a sufficiently
high Nash blowup that they become Cartier. 
This expresses an effective Cartier divisor $-K_U$ as another
Cartier divisor $-D_1$  minus
effective divisors. $D_1$ is the extension of $K_M$ across $U,$ starting
with the smooth locus of $M,$ obtained by pulling back $\Lambda^n\Omega_M$ along $U\to M$.
\m The hypothetical limiting natural number
$-K_U\cdot \gamma$ is the Lefschetz product $-D_1\cdot \gamma$
plus a sum of negative numbers, which must be finite. 
If we insist on thinking in terms of intersection numbers,
we have to allow some of the finitely many components $P$ of $|-D_1|$ containing
our smooth point of  $\gamma(R)$ to contain the whole of it
under all linear equivalences which  keep 
\mbox{$|-K_U|$} effective;
the same integer multiple of each $P\cdot \gamma$ 
occurs on both sides of the equation and it does not matter what value
we should assign to it.
\m
It is much easier to say that we are talking about a very particular finite
ascending chain of locally prinipal ideals defined on finite terms in the tower
of Nash transformations,  which we can pull back along  successive lifts of $\gamma$ to obtain
a finite chain of locally free ideal sheaves. It is just a combinatorial fact
that the series in the theorem is a geometric series when this series for $-K_U$
actually does have finitely many nonzero terms.
\vfill\eject\noindent
The principle which we're applying at this point is just a general priciple.
Beyond considering curves, if we consider any irreducible subvariety $V$ of $M$
containing a smooth point, we can define $V$ to have a `stable formal proper transform'
if the stable locus of proper transforms of $V$ in finite chains of Nash transforms is
proper over $V.$ Then the theorem generalizes to say that once $V$ has a stable formal
proper transform, the necessary and sufficient condition to actually have a proper
transform in $U$ is that the ascending series of divisors $K_0+(K_1-K_0)+(K_2-K_1)+...$
converges in the literal sense once restricted to $V.$
The geometric series involving our relatively baspoint-free divisors
$D_i$ reduces combinatorially to this easy general principle, in other words.

\m In either case, what we have then is an inverse  system of line bundles
which restricts to a finite inverse system, and has as its limit the intersection.
\vfill\eject\noindent Let's represent
the limiting bundle as the restriction of the canonical bundle of $U$
along an actual morphism. 
 A rough argument is that on any bounded open set, the effective Cartier discrepancy
$$K_{i+1}-K_i=D_{(n+2)^{i+1}}-(n+2)D_{(n+2)^i}$$ must restrict
to a principal divisor for sufficiently large $i,$ which can be moved
away from the stable proper transform of $V.$ 
\m This is true, but a more careful consideration shows
that when the divisor restricts to a principal divisor, it actually
only implies that the Nash blowups are finite maps, at points of 
the tansform of $V$ within that bounded neighbourhood. Consideration of the analytic
local rings of Nash blowups of $M$ at such points of the stable proper transform, and the theorem of finitness
of normalization, implies then that each bounded open subset of  the
stable proper transform
is eventually contained in the smooth locus of a finite Nash blowup,
and therefore contained in an open subset of $U,$ which is the union
of these. 
\m This finishes the proof of theorem 1 under the assumption
 $-K_U$ is effective.  Next, one can argue as follows: 
 if we consider any point $p\in M$, there is a neighbourhood of $p$
in $M$ whose smooth Nash manifold has an effective anticanonical divisor.
This can be seen various ways. One way to think of it is that however
singular the point $p$ may be, we can find $n$ flows
which are everywhere defined in a neighbourhood $W$ of $p$ and transverse
 where they cross  
some other point $q$ in that neighbourhood. 
\m For a moment use the letter $U$ to denote the Nash manifold of $W.$
Since $U\to W$ is natural
 the flows lift to flows on $U.$  Evaluating on the wedge product of the corresponding global
vector fields on $U$ gives a map
${\cal O}_U(K_U)\to {\cal O}_U$ showing that  $-K_U$ is effective.
\vfill\eject\noindent   
Let's 
explicate the geometry which allowed this to take place. Thus
let's consider
the particular ${\cal F}_i$ defined earlier which pull back modulo torsion to ${\cal O}(D_i).$  The tuple of vector fields on $W$
gives a map $$\Omega_W/torsion \to B,$$ 
where $B$ is nothing but the section sheaf of a trivial bundle of rank $n.$
The deRham differential gives a connection
$$d:{\cal O}_W\to {\cal O}_W\otimes B$$
of course. 
\m Starting from ${\cal F}={\cal O}_W,$ the sequence of sheaves
$${\cal F}_1, {\cal F}_2, {\cal F}_3,....$$
now each has a corresponding holomorphic connection
$$\nabla:{\cal F}_i\to {\cal F}_i\otimes B$$
and the sheaf of $n$ forms mod torsion in  the $i$'th Nash blowup
$\pi:W'\to W$ is 
$$(\pi^*{\cal F}^{n+1}\pi^*{\cal F}_{(n+2)^i-1}/torsion)^{-1}\pi^*{\cal F}_{(n+2)^i}/torsion.$$
In more detail, let
$$I={\cal F}_1{\cal F}_{(n+2)}...{\cal F}_{(n+2)^{i-1}}.$$
Then
$${\cal F}_{(n+2)^i-1}=I^{n+1},$$
${\cal F}_{(n+2)^{i}}$ is defined to be $\Lambda^{n+1}{\cal P}({\cal F}I)/torsion,$ 
and the $n$ forms on the $i$'th Nash blowup are 
$$(\pi^*{\cal F}\pi^* I/torsion)^{-n-1}\pi^*\Lambda^{n+1}{\cal P}({\cal F}I)/torsion, \ \ \ \ \ \ \ (1)$$
in this case with ${\cal F}={\cal O}_W$
\m
Once  $t_i$ are local sections which span $I$ and $x_i$ are local coordinates  generating ${\cal O}_W$ then
letting the $s_i$ be the products $t_jx_k$ we have that 
$\Lambda^{n+1}{\cal P}(I)/torsion$ is spanned by
the
$$(s_0\oplus \nabla s_0)\wedge...\wedge (s_n\oplus \nabla s_n).$$
The connection $1\oplus \nabla$ is a connection with values in principal
parts rather than differentials, which is essentially a formal object
whose existence is assured because of properties of the pullback of 
a coherent sheaf to its own Fibr\'{e} Vectoriel.
\m This can be expanded in terms of a tangible connection $\nabla $
though as 
$$\sum_{i=0}^n (-1)^is_i\nabla  (s_0)\wedge
...\wedge \widehat{\nabla\ s_i }\wedge...\wedge \nabla s_n \ \ \ \ \ \ \ (2)$$
 $$=s_0...s_n\sum_{i=0}^n (-1)^i\nabla log (s_0)\wedge
...\wedge \widehat{\nabla\ log s_i }
\wedge...\wedge \nabla log s_n 
$$
\noindent
where the hat denotes a deleted term.

\m For any $i,j$ we have
$$\nabla(s_j)=\nabla(s_i{{s_j}\over{s_i}})$$
$$={{s_j}\over{s_i}}\nabla(s_i)+s_i\otimes d({{s_j}\over{s_i}})$$
and so
$$\nabla log \ s_j-\nabla log \ s_i=dlog({{s_j}\over{s_i}}).$$
Using this,  eliminate $\nabla$ from the answer to obtain
an expression which is just one term, not a sum and not involving $\nabla$
$$s_0...s_n\cdot  dlog (s_1/s_0)\wedge...\wedge d log (s_n/s_0).$$
This shows that the answer does not depend on  $\nabla;$ we are
allowed to choose it.
\vfill\eject\noindent When $i$ is not an $(n+2)$ power but has an expansion
$$\sum a_j(n+2)^j$$ with $0\le a_j<(n+2)$   ${\cal F}_i$ is defined to be $\prod {\cal F}_{(n+2)^j}^{a_j}.$ Thus all ${\cal F}_i$ contain ${\cal F}_1^i$ and any
connection $\nabla$ on ${\cal F}_1$ extends to a meromorphic connection
$${\cal F}_i - \buildrel  \nabla \over \to {\cal F}_i\otimes \Omega_W.$$
\m However, the connection which is constructed explicitly 
using a basis of local derivations (which we can even take
to be commuting)  is a holomorphic connection
$${\cal F}_i\to {\cal F}_i\otimes B$$
for all $i.$ The formula (2) is linear of degree $n+1$ with respect to multiplying
the $s_i$ by a meromorphic function, as the later calculation explains. By the product rule the connection
gives $${\cal F}I\to {\cal F}I\otimes B$$ and the coefficients in ${\cal F}I$ which is locally
principal once pulled back to $W'$ pass through 
 (2) and cancel in (1)  leaving
$$\Lambda^n\Omega_{W'}/torsion\to \Lambda^n({\cal O}_{W'}\otimes B)\cong {\cal O}_{W'}.$$
\m 
The transition map

{\normalsize $$(\tau^*({\cal F}{\cal F}_1...{\cal F}_{(n+2)^{i-1}})/torsion)^{-n-1}
\tau^*{\cal F}_{(n+2)^i}/torsion
$$

 \hfil\mbox{$ \ \ \ \ \ \ \ \ \ \ \ \ \ \ \ \ \ \  \ \ \to
(\tau^*({\cal F}{\cal F}_1...{\cal F}_{(n+2)^i})/torsion)^{-n-1}
\tau^*{\cal F}_{(n+2)^{i+1}}/torsion.$}}

\m  for the Nash blowup $\tau:W''\to M$ of $W'$ is the evident rearrangement of terms in

{\normalsize $$({\cal F}{\cal F}_1...{\cal F}_{(n+2)^i})^{n+1}{\cal F}_{(n+2)^i}
\to ({\cal F}{\cal F}_1...{\cal F}_{(n+2)^{i-1}})^{n+1}{\cal F}_{(n+2)^{i+1}}.$$
}

\noindent induced by carrying 
${\cal F}_{(n+2)^i}^{n+2}\to {\cal F}_{(n+2)^{i+1}}.$

\vfill\eject\noindent
The inclusion commutes with this
transition map;
and in the limit on $U$ it embeds
${\cal O}_U(K_U)\to {\cal O}_U$ 
as the coherent sheaf
of ideals defining
 an effective anticanonical divsor.
Thus if one doesn't wish to trust an argument based on
symmetry alone, it is possible to prove that the anticanonical
divisor of the Nash manifold is linearly equivalent to a particular effective 
divisor when restricted to the inverse
image of sufficiently small open neighbourhoods of a point.  QED
\m{\bf Remark.}  It is probably nicer to show
 that there is  always a locally principal  sheaf $J$ on
$M$ conducting $n$ forms into an ideal sheaf on $U,$ analogous
to how $dlog(s_1/s_0)$ has no zeroes on the projective line being a section
of ${\cal O}(-[s_1]-[s_0]). $
Maybe ${\cal F}\cong {\cal O}_W$ on open
sets $W$ are restrictions ${\cal F}(W)$ of some  ${\cal F}$  with
 connection on $M.$ $J$ might be $\Lambda^{-n}B$ for $B$ containing
the one forms of $M,$ with the meromorphic extension 
${\cal F}_i\to {\cal F}_i\otimes B$ holomorphic everywhere. 
\vfill\eject\noindent
References.

\noindent
1. Ian Hislop, The Olden Days, Green imagined land, Episode 3, BBC 2, 23 April, (2014)

\noindent
2. B. Mandelbrot Mandelbrot,  The Fractal Geometry of Nature. W.H. Freeman and Co (1983)

\noindent
3. D. Mumford, A.  Grothendieck, C.-L. Chai, Selected papers on the classification of varieties and moduli spaces
 Springer (2004)

\noindent
4. E. Witten, Physics and Geometry, ICM (1986)

\noindent 
5 The minimal model program (circa 1979)

\noindent
6. The genetic literacy project, The more biotech science you know, the less you fear GMO crops, study finds
(see also McPhetres, Journal of Environmental Psychology,  May 29 (2019)

\noindent
7. H. Hironaka,   Resolution of Singularities of an Algebraic Variety Over a Field of Characteristic Zero: I, II, Annals of Mathematics (1964)

\noindent
8. P. Perez, B. Teissier, Toric geometry and semple-nash modification,
Revista de la Real Academia de Ciencias (2014)

\noindent
9.  J. G. Semple, Some investigations in the geometry of curve and surface elements, Proc. London Math. Soc.
 (1954) 

\noindent
10. M. Spivakovski, conversation at conference in Grenoble organized by
V. Cossart 

\noindent
11. S. Abyankhar (Ibid.)

\noindent
12. S.  Abyankhar, conversation at conference at Warwick

\noindent
13. G. Gonzalez-Sprinberg, Eventails en dimension deux et transforme de Nash,  Publications
de l’Ecole Normale Superieure  (1977)

\noindent
14. A. Atanasov, C. Lopez, A. Perry, N. Proudfoot, M. Thaddeus,
Resolving Toric Varieties with Nash Blowups, Experimental Mathematics (2011)

\noindent
15. S. Encinas, O. Villamayor, Good points and constructive resolution of singularities. Acta
Math. 181 (1998) 

\noindent
16. J. Wlodzarcek  Simple Hironaka resolution in characteristic zero, J. Amer. Math.
Soc. 18 (2005),  

\noindent
17. J. Kollar, Lectures on resolution of singularities, Annals of Mathematics studies (2007)

\noindent
18. M. McQuillan Very fast, very functorial, and very easy resolution of singularities, arXiv (2019)

\noindent
19. D. AAbramovich, M. Temkin, J. Wlodarczyk, Functorial resolution
by weighted blow-ups (2019)  

\noindent
20. J. P. Serre, Algebre Local, Multiplicites, Springer Lecture notes 11

\noindent
21.A.  Dold, Homology of symmetric products and other functors of complexes. Ann. of
Math. (1958)

\noindent
22. D. Gabriel, A. Zisman, Calculus of fractions and homotopy theory, Ergebnisse der Mathematik und ihrer Grenzgebiete (1967)

\noindent
23. J. Moody, On resolving Singularities LMS (2001)

\noindent
24. J. Moody, Finite generation and the Gauss process (preprint)

\noindent
25. A natural Chern character, J. Moody (preprint)

\noindent
26. Functorial affinization of Nash's manifold (preprint)

\end{document}